\newif\ifpagetitre            \pagetitretrue
\newtoks\hautpagetitre        \hautpagetitre={\hfil}
\newtoks\baspagetitre         \baspagetitre={\hfil}
\newtoks\auteurcourant        \auteurcourant={\hfil}
\newtoks\titrecourant         \titrecourant={\hfil}

\newtoks\hautpagegauche       \newtoks\hautpagedroite
\hautpagegauche={\hfil\the\auteurcourant\hfil}
\hautpagedroite={\hfil\the\titrecourant\hfil}

\newtoks\baspagegauche \baspagegauche={\hfil\tenrm\folio\hfil}
\newtoks\baspagedroite \baspagedroite={\hfil\tenrm\folio\hfil}

\headline={\ifpagetitre\the\hautpagetitre
\else\ifodd\pageno\the\hautpagedroite
\else\the\hautpagegauche\fi\fi}

\footline={\ifpagetitre\the\baspagetitre
\global\pagetitrefalse
\else\ifodd\pageno\the\baspagedroite
\else\the\baspagegauche\fi\fi}

\vsize=9.0in\voffset=1cm
%\nopagenumbers
\looseness=2

% formato JGR

\message{fonts,}

\font\tenrm=cmr10
\font\ninerm=cmr9
\font\eightrm=cmr8
\font\teni=cmmi10
\font\ninei=cmmi9
\font\eighti=cmmi8
\font\ninesy=cmsy9
\font\tensy=cmsy10
\font\eightsy=cmsy8
\font\tenbf=cmbx10
\font\ninebf=cmbx9
\font\tentt=cmtt10
\font\ninett=cmtt9

\font\ninesl=cmsl9
\font\eightsl=cmsl8

\font\nineit=cmti9
\font\eightit=cmti8
 % for titles

\skewchar\ninei='177 \skewchar\eighti='177
\skewchar\ninesy='60 \skewchar\eightsy='60

\def\eightpoint{\def\rm{\fam0\eightrm} % pasa a font de 8 puntos
\normalbaselineskip=9pt
\normallineskiplimit=-1pt
\normallineskip=0pt

\textfont0=\eightrm \scriptfont0=\sevenrm \scriptscriptfont0=\fiverm
\textfont1=\ninei \scriptfont1=\seveni \scriptscriptfont1=\fivei
\textfont2=\ninesy \scriptfont2=\sevensy \scriptscriptfont2=\fivesy
\textfont3=\tenex \scriptfont3=\tenex \scriptscriptfont3=\tenex
\textfont\itfam=\eightit  \def\it{\fam\itfam\eightit} % \it is family 4
\textfont\slfam=\eightsl \def\sl{\fam\slfam\eightsl} % \sl is family 5

\setbox\strutbox=\hbox{\vrule height6pt depth2pt width0pt}%
\normalbaselines \rm}

\def\ninepoint{\def\rm{\fam0\ninerm} % pasa a font de 9 puntos
\textfont0=\ninerm \scriptfont0=\sevenrm \scriptscriptfont0=\fiverm
\textfont1=\ninei \scriptfont1=\seveni \scriptscriptfont1=\fivei
\textfont2=\ninesy \scriptfont2=\sevensy \scriptscriptfont2=\fivesy
\textfont3=\tenex \scriptfont3=\tenex \scriptscriptfont3=\tenex
\textfont\itfam=\nineit  \def\it{\fam\itfam\nineit} % \it is family 4
\textfont\slfam=\ninesl \def\sl{\fam\slfam\ninesl} % \sl is family 5
\textfont\bffam=\ninebf \scriptfont\bffam=\sevenbf
\scriptscriptfont\bffam=\fivebf \def\bf{\fam\bffam\ninebf} % \bf is family 6
\textfont\ttfam=\ninett \def\tt{\fam\ttfam\ninett} % \tt is family 7

\normalbaselineskip=11pt
\setbox\strutbox=\hbox{\vrule height8pt depth3pt width0pt}%
\let \smc=\sevenrm \let\big=\ninebig \normalbaselines
\parindent=1em
\rm}

\def\tenpoint{\def\rm{\fam0\tenrm} % font de 10 points
\textfont0=\tenrm \scriptfont0=\ninerm \scriptscriptfont0=\fiverm
\textfont1=\teni \scriptfont1=\seveni \scriptscriptfont1=\fivei
\textfont2=\tensy \scriptfont2=\sevensy \scriptscriptfont2=\fivesy
\textfont3=\tenex \scriptfont3=\tenex \scriptscriptfont3=\tenex
\textfont\itfam=\nineit  \def\it{\fam\itfam\nineit} % \it is family 4
\textfont\slfam=\ninesl \def\sl{\fam\slfam\ninesl} % \sl is family 5
\textfont\bffam=\ninebf \scriptfont\bffam=\sevenbf
\scriptscriptfont\bffam=\fivebf \def\bf{\fam\bffam\tenbf} % \bf is family 6
\textfont\ttfam=\tentt \def\tt{\fam\ttfam\tentt} % \tt is family 7

\normalbaselineskip=11pt
\setbox\strutbox=\hbox{\vrule height8pt depth3pt width0pt}%
\let \smc=\sevenrm \let\big=\ninebig \normalbaselines
\parindent=1em
\rm}

\message{fin format jgr}
%end format jgr

\hautpagegauche={\hfill\ninerm\the\auteurcourant}
\hautpagedroite={\ninerm\the\titrecourant\hfill}
\auteurcourant={R.G. Novikov}
\titrecourant={On stable determination of potential by boundary measurements}

\magnification=1200
\font\Bbb=msbm10
\def\R{\hbox{\Bbb R}}
\def\C{\hbox{\Bbb C}}
\def\N{\hbox{\Bbb N}}
\def\S{\hbox{\Bbb S}}
\def\pa{\partial}

\vskip 2 mm
\centerline{\bf On stable determination of potential by boundary
measurements}

\vskip 2 mm
\centerline{\bf R. G. Novikov}
\vskip 4 mm

\noindent
{\ninerm CNRS, Laboratoire de Math\'ematiques Jean Leray (UMR 6629),
Universit\'e de Nantes,}

\noindent
{\ninerm BP 92208, F-44322, Nantes cedex 03, France}

\noindent
{\ninerm e-mail: novikov@math.univ-nantes.fr}

\vskip 2 mm
{\bf Abstract.}
We give new stability estimates for the Gel'fand-Calderon inverse
boundary value problem.

\vskip 2 mm
{\bf 1. Introduction}

Consider the equation
$$-\Delta\psi + v(x)\psi=0,\ \ x\in D, \eqno(1.1)$$
where
$$\eqalign{
&D\ \ {\rm is\ an\ open\ bounded\ domain\ in}\ \ \R^d,\cr
&\pa D\in C^2,\ \ v\in L^{\infty}(D),\ \ d\ge 2.\cr}\eqno(1.2)$$
We assume also that
$$\eqalign{
&0\ \ {\rm is\ not\ a\ Dirichlet\ eigenvalue\ for}\cr
&{\rm the\ operator}\ \ -\Delta+v\ \ {\rm in}\ \ D.\cr}\eqno(1.3)$$

Equation (1.1) arises, in particular, in quantum mechanics, acoustics,
electrodynamics. Formally, (1.1) looks as the Schr\"odinger equation with
potential $v$ at zero energy.

Consider the map $\Phi$ such that
$${\pa\psi\over \pa\nu}\big|_{\pa D}=\Phi(\psi\big|_{\pa D}) \eqno(1.4)$$
for all sufficiently regular solutions $\psi$ of (1.1) in
$\bar D=D\cup\pa D$, where $\nu$ is the outward normal to $\pa D$. The map
$\Phi$ is called the Dirichlet-to-Neumann map for equation (1.1) and is
considered as boundary measurements for (1.1).

We consider the following inverse boundary value problem for equation (1.1):

\vskip 2 mm
{\bf Problem 1.1.}
Given $\Phi$, find $v$.

This problem can be considered as the Gel'fand inverse boundary value problem
 for the Schr\"odinger equation at zero energy (see [G], [No1]). This
problem can be also considered as a generalization of the Calderon problem
of the electrical impedance tomography (see [C], [SU], [No1]).
Concerning results given in the literature on Problem 1.1 (in its Calderon
or Gel'fand form ) see  [KV], [SU], [HN] (note added in proof),
[No1], [Al], [Na1], [Na2], [BU], [P], [Ma], [No2], [No4], [Am]   and
references therein.

In the present article we show that the Alessandrini stability estimates
of [Al] for Problem 1.1 in dimension  $d\ge 3$ (see Theorem 2.1 of the
next section) admit some principle improvement. Our new stability estimates
(see Theorem 2.2 of the next section) are obtained by methods developed in
[No2], [No3], [No4]. These methods include, in particular: (1) the
$\bar\pa$- approach to inverse "scattering" at zero energy in dimension
$d\ge 3$, going back to [BC], [HN], and (2) the reduction of Problem 1.1
to inverse "scattering" at zero energy, going back to [No1].

The present article is organized as follows.
In Section 2, we formulate and discuss old and new stability estimate for
Problem 1.1. In Section 3, we remind (a) definition and some properties
of the Faddeev functions and (b) formulation of the inverse "scattering"
problem  for the Schr\"odinger equation at zero energy (Problem 3.1).
In Section 4, we remind formulas and equations of [No1], [No2] reducing
Problem 1.1 to Problem 3.1. In Section 5, we remind an approximate
reconstruction method of [No1] for Problem 1.1. In Section 6 we prove
Theorem 2.2 in the Born approximation.

\vskip 2 mm
{\bf 2. Stability estimates}

We assume for simplicity that
$$\eqalign{
&D\ \ {\rm is\ an\ open\ bounded\ domain\ in}\ \ \R^d,\ \pa D\in C^2,\cr
&v\in W^{m,1}(\R^d)\ \ {\rm for\ some}\ m>d,\ supp\,v\subset D,\ d\ge 2,\cr}
\eqno(2.1)$$
where
$$W^{m,1}(\R^d)=\{v:\ \pa^Jv\in L^1(\R^d),\ |J|\le m\},\ \ m\in\N\cup 0,
\eqno(2.2)$$
where
$$J\in (\N\cup 0)^d,\ \ |J|=\sum_{i=1}^dJ_i,\ \
\pa^Jv(x)={\pa^{|J|}v(x)\over \pa x_1^{J_1}\ldots \pa x_d^{J_d}}.$$
Let
$$\|v\|_{m,1}=\max_{|J|\le m}\|\pa^Jv\|_{L^1(\R^d)}.\eqno(2.3)$$
Let
$$\eqalign{
&\|A\|\ \ {\rm denote\ the\ norm\ of\ an\ operator}\cr
&A:\ L^{\infty}(\pa D)\to L^{\infty}(\pa D).\cr} \eqno(2.4)$$
We remind that if $v_1$, $v_2$ are potentials satisfying (1.2), (1.3),
where $D$ is fixed, then
$$\Phi_1-\Phi_2\ \ {\rm is\ a\ compact\ operator\ in}\ \ L^{\infty}(\pa D),
\eqno(2.5)$$
where $\Phi_1$, $\Phi_2$ are the DtN maps for $v_1$, $v_2$ respectively,
see [No1], [No2]. Note also that $(2.1)\Rightarrow (1.2)$.

\vskip 2 mm
{\bf Theorem 2.1} (variation of the result of [Al]).
{\it  Let conditions} (1.3), (2.1) {\it hold for potentials} $v_1$
{\it and} $v_2$, {\it where} $D$ {\it is fixed}, $d\ge 3$. {\it Let}
$\|v_j\|_{m,1}\le R$, $j=1,2$, {\it for some} $R>0$. {\it Let} $\Phi_1$,
$\Phi_2$ {\it denote the DtN maps for} $v_1$, $v_2$, {\it respectively. Then}
$$\|v_1-v_2\|_{L^{\infty}(D)}\le
C_1(\ln(1+\|\Phi_1-\Phi_2\|^{-1}))^{-\alpha_1},\eqno(2.6)$$
{\it where} $C_1=C_1(R,D,m)$, $\alpha_1=(m-d)/m$, $\|\Phi_1-\Phi_2\|$
{\it is defined according to} (2.4).

Theorem 2.1 follows from formulas (3.9)-(3.11), (4.1) (of Sections 3 and 4).

A disadvantage of estimate (2.6) is that
$$\alpha_1<1\ \ {\rm for\ any}\ \  m>d\ \ {\rm even\ if}\ \ m\ \ {\rm is\
 very\ great}.\eqno(2.7)$$

\vskip 2 mm
{\bf Theorem 2.2.}
{\it Let the assumptions of Theorem 2.1 hold. Then}
$$\|v_1-v_2\|_{L^{\infty}(D)}\le
C_2(\ln(1+\|\Phi_1-\Phi_2\|^{-1}))^{-\alpha_2},\eqno(2.8)$$
{\it where} $C_2=C_2(R,D,m)$, $\alpha_2=m-d$, $\|\Phi_1-\Phi_2\|$
{\it is defined according to} (2.4).

A principal advantage of estimate (2.8) in comparison with (2.6) is that
$$\alpha_2\to +\infty\ \ {\rm as}\ \  m\to +\infty, \eqno(2.9)$$
in contrast with (2.7).

In the Born approximation, that is in the linear approximation near zero
potential, Theorem 2.2 is proved in Section 6.

For sufficiently small $R$ in dimension $d=3$, Theorem 2.2 follows
from (3.9) (of Section 3) and results of [No2], [No4].
The scheme of our proof for this case is, actually, similar to the scheme
of our proof for the case of the Born approximation. The main difference
is that instead of the inverse Fourier transform (used in Section 6) we use
now the zero-energy inverse "backscattering" transform of [No4]. We plan to
give this "nonlinear" "small-norm" proof in a separate article.

In the general case, the proof of Theorem 2.2 is not completed yet. However,
except restrictions in time,  we see no difficulties for
completing this proof by methods of [No2], [No3], [No4].

\vskip 2 mm
{\bf 3. Faddeev functions}

We consider the Faddeev functions $G$, $\psi$ and $h$ (see [F1], [F2],
[HN], [No2]):
$$\eqalignno{
&G(x,k)=e^{ikx}g(x,k),\ \
g(x,k)=-\bigl({1\over 2\pi}\bigr)^d\int_{\R^d}{e^{i\xi x}d\xi\over
{\xi^2+2k\xi}},&(3.1)\cr
&\psi(x,k)=e^{ikx}+\int_{\R^d}G(x-y,k)v(y)\psi(y,k)dy,&(3.2)\cr}$$
where $x\in\R^d$, $k\in\Sigma$,
$$\eqalignno{
&\Sigma=\{k\in\C^d,\ \ k^2=k_1^2+\ldots +k_d^2=0\};&(3.3)\cr
&h(k,l)=\bigl({1\over 2\pi}\bigr)^d\int_{\R^d}e^{-ilx}v(x)\psi(x,k)dx,&(3.4)
\cr}$$
where $(k,l)\in\Theta$,
$$\Theta=\{k\in\Sigma,\ \ l\in\Sigma:\ Im\,k=Im\,l\}.\eqno(3.5)$$
We remind that:
$$\Delta G(x,k)=\delta(x),\ \ x\in\R^d,\ \ k\in\Sigma;\eqno(3.6)$$
formula (3.2) at fixed $k$ is considered as an equation for
$$\psi=e^{ikx}\mu(x,k),\eqno(3.7)$$
where $\mu$ is sought in $L^{\infty}(\R^d)$; as a corollary of (3.2), $\psi$
satisfies (1.1); $h$ of (3.4) is a generalized "scattering" amplitude
in the complex domain at zero energy.

Note that, actually, $G$, $\psi$, $h$ of (3.1)-(3.5) are zero energy
restrictions of functions introduced by Faddeev as extentions to the
complex domain of some functions of the classical scattering theory
for the Schr\"odinger equation at positive energies. In addition,
$G$, $\psi$, $h$ in their zero energy restriction were considered for the
first time in [BC]. The Faddeev functions $G$, $\psi$, $h$ were, actually,
rediscovered in [BC].

We remind also that, under the assumptions of Theorem 2.1
$$\mu(x,k)\to 1\ \ {\rm as}\ \ |Im\,k|\to\infty\ \ ({\rm uniformly\ in}\ \ x)
\eqno(3.8)$$
and, for any $c>1$,
$$|\mu(x,k)|< c\ \ {\rm for}\ \ |Im\,k|\ge\rho_1(R,D,m,c),\eqno(3.9)$$
where $x\in\R^d$, $k\in\Sigma$;
$$\hat v(p)=\lim_{\scriptstyle (k,l)\in\Theta,\ k-l=p\atop\scriptstyle
|Im\,k|=|Im\,l|\to\infty}h(k,l)\ \ {\rm for\ any}\ \ p\in\R^d,\eqno(3.10)$$
$$\eqalign{
&|\hat v(p)-h(k,l)|\le {C_3(D,m)R^2\over \rho}\ \ {\rm for}\ \
(k,l)\in\Theta,\ p=k-l,\cr
&|Im\,k|=|Im\,l|=\rho\ge\rho_2(R,D,m),\cr}\eqno(3.11)$$
where
$$\hat v(p)=\bigl({1\over 2\pi}\bigr)^d\int_{\R^d}e^{ipx}v(x)dx,\ \ p\in\R^d.
\eqno(3.12)$$

Results of the type (3.8), (3.9) go back to [BC]. Results of the type
(3.10), (3.11) (with less precise right-hand side in (3.11)) go back to [HN].
Estimates (3.8), (3.11) are  related also with some important
$L_2$-estimate going back to [SU] on the Green function $g$ of (3.1).

For more information on properties of the Faddeev functions
$G$, $\psi$, $h$, see [HN], [No2], [No4] and references therein.

In the next section we remind that Problem 1.1 (of Introduction) admits a
reduction to the following inverse "scattering" problem:

\vskip 2 mm
{\bf Problem 3.1.}
Given $h$ on $\Theta$, find $v$ on $\R^d$.

\vskip 2 mm
{\bf 4. Reduction of [No1], [No2]}

Let conditions (1.2), (1.3) hold for potentials $v_1$ and $v_2$, where
$D$ is fixed. Let $\Phi_i$, $\psi_i$, $h_i$ denote the DtN map $\Phi$
and the Faddeev functions $\psi$, $h$ for $v=v_i$, $i=1,2$. Let also
$\Phi_i(x,y)$ denote the Schwartz kernel $\Phi(x,y)$ of the integral operator
$\Phi$ for $v=v_i$, $i=1,2$. Then (see [No2] for details):
$$h_2(k,l)-h_1(k,l)=\bigl({1\over 2\pi}\bigr)^d
\int\limits_{\pa D}\!\!\int\limits_{\pa D}
\psi_1(x,-l)(\Phi_2-\Phi_1)(x,y)\psi_2(y,k)dydx,\eqno(4.1)$$
where $(k,l)\in\Theta$;
$$\eqalignno{
&\psi_2(x,k)=\psi_1(x,k)+\int_{\pa D}A(x,y,k)\psi_2(y,k)dy,\ x\in\pa D,
&(4.2a)\cr
&A(x,y,k)=\int_{\pa D}R_1(x,z,k)(\Phi_2-\Phi_1)(z,y)dz,\ x,y\in\pa D,
&(4.2b)\cr
&R_1(x,y,k)=G(x-y,k)+\int_{\R^d}G(x-z,k)v_1(z)
R_1(z,y,k)dz,\ x,y\in\R^d,
&(4.3)\cr}$$
where $k\in\Sigma$. Note that: (4.1) is an explicit formula, (4.2a) is
considered as an equation for finding $\psi_2$ on $\pa D$ from $\psi_1$ on
$\pa D$ and $A$ on $\pa D\times\pa D$ for each fixed $k$, (4.2b)
is an explicit formula, (4.3) is an equation for finding $R_1$ from $G$
and $v_1$, where $G$ is the function of (3.1).

Note that formulas and equations (4.1)-(4.3) for $v_1\equiv 0$ were given
in [No1] (see also [HN] (Note added in proof), [Na1], [Na2]).
In this case $h_1\equiv 0$, $\psi_1=e^{ikx}$, $R_1=G(x-y,k)$. Formulas and
equations (4.1)-(4.3) for the general case were given in [No2].

Formulas and equations (4.1)-(4.3) with fixed background potential $v_1$
reduce Problem 1.1 (of Introduction) to Problem 3.1 (of Section 3).

\vskip 2 mm
{\bf 5. Reconstruction of [No1] in the Born approximation}

In the Born approximation, that is in the linear approximation near zero
potential, we have that
$$\eqalignno{
&h(k,l)\approx\hat v(k-l),&(5.1)\cr
&h(k,l)\approx\bigl({1\over 2\pi}\bigr)^d
\int\limits_{\pa D}\!\!\int\limits_{\pa D}
e^{-ilx}(\Phi-\Phi_0)(x,y)e^{iky}dxdy,&(5.2)\cr}$$
where $(k,l)\in\Theta$, $\hat v$ is defined by (3.12), $\Phi_0$ denotes
the DtN map for $v\equiv 0$.

Formulas (5.1), (5.2) follow from (3.1)-(3.4) and (4.1). Formulas (5.1), (5.2)
 imply, in particular, that
$$\eqalignno{
&\hat v(p)\approx\bigl({1\over 2\pi}\bigr)^d
\int\limits_{\pa D}\!\!\int\limits_{\pa D}
e^{-il(p)x}(\Phi-\Phi_0)(x,y)e^{ik(p)y}dxdy,&(5.3)\cr
&k(p)={p\over 2}+i{|p|\over 2}\gamma(p),\ \
l(p)=-{p\over 2}+i{|p|\over 2}\gamma(p),\ \ p\in\R^d,&(5.4a)\cr}$$
where $\gamma(p)$ is a piecewise continuous function of
$p\in\R^d$ with values in $\S^{d-1}$ and such that
$$\gamma(p)p=0,\ \ p\in\R^d.\eqno(5.4b)$$
One can see that formula (5.3) gives a reconstruction method for Problem 1.1,
$d\ge 2$, in the Born approximation.

An approximate reconstruction method based on (5.1), (5.2) for Problem 1.1
in dimension $d\ge 2$ was proposed for the first time in [No1].

In the next section we show that, in the Born approximation, Theorem 2.2
(of Section 2) follows, actually, from (5.3).

\vskip 2 mm
{\bf 6. Proof of Theorem 2.2 in the Born approximation}

We have that
$$\eqalignno{
&v_1(x)-v_2(x)=\biggl(\int\limits_{|p|<\rho}+\int\limits_{|p|>\rho}\biggr)
e^{-ipx}(\hat v_1(p)-\hat v_2(p))dp,&(6.1)\cr
&|v_1(x)-v_2(x)|\le I_1(\rho)+I_2(\rho),&(6.2a)\cr
&I_1(\rho)=\int\limits_{|p|<\rho}|\hat v_1(p)-\hat v_2(p)|dp,&(6.2b)\cr
&I_2(\rho)=\int\limits_{|p|>\rho}|\hat v_1(p)-\hat v_2(p)|dp,&(6.2c)\cr}$$
where $x\in\R^d$, $\rho>0$.

The assumptions $\|v_j\|_{m,1}\le R$, $j=1,2$, imply that
$$|\hat v_1(p)-\hat v_2(p)|\le {C_4(d,m)R\over (1+|p|)^m},\ \ p\in\R^d.
\eqno(6.3)$$
Using (5.3) we obtain that
$$\eqalign{
&|\hat v_1(p)-\hat v_2(p)|\le C_5(D)e^{L\rho}\|\Phi_1-\Phi_2\|,\ \
|p|\le\rho,\cr
&C_5(D)=\bigl({1\over 2\pi}\bigr)^d\int\limits_{\pa D}dx,\ \
L=\max_{x\in D}|x|,\cr}\eqno(6.4)$$
where $\|\cdot\|$ is defined according to (2.4).

Formulas (6.2b), (6.4) imply that
$$\eqalign{
&I_1(\rho)\le C_6(D)\rho^de^{L\rho}\|\Phi_1-\Phi_2\|\le
C_6(D)e^{L_1\rho}\|\Phi_1-\Phi_2\|,\ \ \rho>0,\cr
&C_6(D)=C_5(D)\int\limits_{\theta\in\S^{d-1}}d\theta,\ \
L_1=L+d.\cr}\eqno(6.5)$$
Formulas (6.2c), (6.3) imply that
$$I_2(\rho)\le C_7(d,m)R\rho^{-(m-d)},\ \ \rho>0.\eqno(6.6)$$
Let
$$\eqalign{
&\alpha\in ]0,1[\ \ {\rm be\ fixed},\cr
&\|\Phi_1-\Phi_2\|=\delta,\ \ \rho=\lambda\ln(1+\delta^{-1}),\ \
\lambda={{1-\alpha}\over L_1}.\cr}\eqno(6.7)$$
Then, due to (6.5), (6.6),
$$\eqalign{
&I_1(\rho)\le C_6\delta(1+\delta^{-1})^{\lambda L_1}=
C_6(1+\delta)^{1-\alpha}\delta^{\alpha},\cr
&I_2(\rho)\le C_7R(\lambda\ln(1+\delta^{-1}))^{-(m-d)}.\cr}\eqno(6.8)$$
Estimate (2.8) for $\delta=\|\Phi_1-\Phi_2\|\le 1/2$ follows from (6.2a),
(6.8). Estimate (2.8) in the general case (with modified $C_2$)
follows from (2.8) for $\delta\le 1/2$ and the assumptions that
$\|v_j\|_{L^1(\R^d)}\le R$, $j=1,2$.

Thus, in the  Born approximation Theorem 2.2 is proved.
(This proof is valid also for $d=2$).

\vskip 4 mm
{\bf References}
\vskip 2 mm
\item{[  Al]} G.Alessandrini, {\it Stable determination of conductivity
by boundary measurements}, Appl. Anal. {\bf 27} (1988), 153-172.
\item{[  Am]} H.Ammari, {\it An Introduction to Mathematics of Emerging
Biomedical Imaging},

\item{      } Springer, Berlin, 2007
\item{[  BC]} R.Beals and R.R.Coifman, {\it Multidimensional inverse
scattering and nonlinear partial differential equations}, Proc. Symp. Pure
Math. {\bf 43} (1985), 45-70.
\item{[  BU]} R.M.Brown and G.Uhlmann, {\it Uniqueness in the inverse
conductivity problem for nonsmooth conductivities in two dimensions},
Comm. Partial Diff. Eq. {\bf 22} (1997), 1009-1027.
\item{[   C]} A.-P.Calder\'on, {\it On an inverse boundary value problem},
Seminar on Numerical Analysis and its Applications to Continuum Physics
(Rio de Janeiro, 1980), pp.65-73, Soc. Brasil. Mat. Rio de Janeiro, 1980.
\item{[  F1]} L.D.Faddeev, {\it Growing solutions of the Schr\"odinger
equation}, Dokl. Akad. Nauk SSSR {\bf 165} (1965), 514-517 (in Russian);
English Transl.: Sov. Phys. Dokl. {\bf 10} (1966), 1033-1035.
\item{[  F2]} L.D.Faddeev, {\it Inverse problem of quantum scattering theory
II}, Itogi Nauki i Tekhniki, Sovr. Prob. Math. {\bf 3} (1974), 93-180
(in Russian); English Transl.: J.Sov. Math. {\bf 5} (1976), 334-396.
\item{[   G]} I.M.Gelfand, {\it Some problems of functional analysis and
algebra}, Proceedings of the International Congress of Mathematicians,
Amsterdam, 1954, pp.253-276.
\item{[  HN]} G.M.Henkin and R.G.Novikov, {\it The $\bar\pa$- equation in the
multidimensional inverse scattering problem}, Uspekhi Mat. Nauk {\bf 42(3)}
(1987), 93-152 (in Russian); English Transl.: Russ. Math. Surv. {\bf 42(3)}
(1987), 109-180.
\item{[  KV]} R.Kohn and M.Vogelius, {\it Determining conductivity by
boundary measurements II, Interior results}, Comm. Pure Appl. Math. {\bf 38}
(1985), 643-667.
\item{[  Ma]} N.Mandache, {\it Exponential instability in an inverse problem
for the Schr\"odinger equation}, Inverse Problems {\bf 17} (2001), 1435-1444.
\item{[ Na1]} A.I.Nachman, {\it Reconstructions from boundary measurements},
Ann. Math. {\bf 128} (1988), 531-576.
\item{[ Na2]} A.I.Nachman, {\it Global uniqueness for a two-dimensional
inverse boundary value problem}, Ann, Math. {\bf 142} (1995), 71-96.
\item{[ No1]} R.G.Novikov, {\it Multidimensional inverse spectral problem
for the equation $-\Delta\psi+(v(x)-Eu(x))\psi=0$}, Funkt. Anal. i Pril.
{\bf 22(4)} (1988), 11-22 (in Russian); English Transl.: Funct. Anal. and
Appl. {\bf 22} (1988), 263-272.
\item{[ No2]} R.G.Novikov, {\it Formulae and equations for finding scattering
data from the Dirichlet-to-Neumann map with nonzero background potential},
Inverse Problems {\bf 21} (2005), 257-270.
\item{[ No3]} R.G.Novikov, {\it The $\bar\pa$- approach to approximate inverse
scattering at fixed energy in three dimensions}, International Mathematics
Research Papers, {\bf 2005:6}, (2005), 287-349.
\item{[ No4]} R.G.Novikov, {\it On non-overdetermined inverse scattering
at zero energy in three dimensions}, Ann. Scuola Norm. Sup. Pisa Cl. Sci.
{\bf 5} (2006), 279-328
\item{[   P]} V.P.Palamodov, {\it Gabor analysis of the continuum model
for impedance tomography}, Ark.Mat. {\bf 40} (2002), 169-187
\item{[  SU]} J.Sylvester and G.Uhlmann, {\it A global uniqueness theorem
for an inverse boundary value problem}, Ann. Math. {\bf 125} (1987),
153-169

\end